


\parindent = 0 pt
\baselineskip = 16 pt
\parskip = \the\baselineskip

\font\AMSBoldBlackboard = msbm10

\def\RR{{\hbox{\AMSBoldBlackboard R}}}
\def\CC{{\hbox{\AMSBoldBlackboard C}}}
\def\AA{{\hbox{\AMSBoldBlackboard A}}}
\def\QQ{{\hbox{\AMSBoldBlackboard Q}}}

\def\ZZ{{\hbox{\AMSBoldBlackboard Z}}}

\settabs 12\columns
\rightline{13-FEB-1999}
\rightline{math.NT/9902080}
\vskip 1 true in
{\bf \centerline{THE EXPLICIT FORMULA AND THE CONDUCTOR OPERATOR}}
\vskip 0.5 true in
\centerline{Jean-Fran\c{c}ois Burnol}
\par
\centerline{February 1999}
\par
I give a new derivation of the Explicit Formula for an arbitrary number field and abelian Dirichlet-Hecke character, which treats all primes in exactly the same way, whether they are discrete or archimedean, and also ramified or not. This is followed with a local study of a Hilbert space operator, the ``conductor operator'', which is expressed as $H = \log(|x|) + \log(|y|)$ (where $x$ and $y$ are Fourier dual variables on a $\nu$-adic completion of the number field). I also study the commutator operator $K = i\ [\log(|y|),\ \log(|x|)]$ (which shares with $H$ the property of complete dilation invariance, and turns out to be bounded), as well as the higher commutator operators. The generalized eigenvalues of these operators are given by the derivatives on the critical line of the Tate-Gel'fand-Graev Gamma function, which itself is in fact closely related to the additive Fourier Transform viewed in multiplicative terms. This spectral analysis is thus a natural continuation to Tate's Thesis in its local aspects.\par

\vfill
{\parskip = 0 pt
62 rue Albert Joly\par
F-78000 Versailles\par
France\par}

\eject

{\bf TABLE OF CONTENTS\par
\par
A. Introduction\par
B. Necessary facts from Tate's Thesis\par
C. The abelian Explicit Formula\par
D. Fourier transform of $-\log (|y|)$\par
E. The Explicit Formula and the conductor operator\par
F. Concrete spectral analysis on $\QQ_p$\par
G. Transformation Theory and an intrinsic spectral analysis\par
H. More on homogeneous distributions and Gamma functions\par
I. Taking derivatives and getting conductor and commutators\par
J. More on the Hilbert Space story\par
K. Notes\par
References\par}
\vfill \eject

{\bf A. Introduction}\par

Ever since Tate's Thesis [Ta50] (see also [Iw52]) there has always been a tension between the additive and multiplicative structures in the adelic approach to the Riemann Zeta function and its (abelian) generalizations. Tate gave a proof of the functional equation which reduced it to the computation of the additive Fourier transform of a multiplicatively homogeneous distribution.

In a paper dating back to about the same period Weil [We52] puts forward an identity $$Z(K, \chi )(g) = W(K, \chi )(g)$$ which encompasses most abelian ``Explicit Formulae''. Here $K$ is a number field, $\chi$ a unitary Dirichlet-Hecke character on its idele group and $g$ a function on the positive real line, which we will assume to be smooth and compactly supported (although Weil considered functions of a much wider class). The left-hand-side $Z(K, \chi )(g)$ is obtained by adding the values of the Mellin transform
$$\widehat{g}(s) = \int_{0}^{\infty} g(u)u^s \, { du \over u}$$
of $g$ at the poles (counted with positive multiplicity) and the zeros (counted with negative multiplicities) of the $L$-function $L(K, \chi, s)$. The right-hand-side $W(K, \chi)$ has an additive decomposition  over the prime spots of $K$. Weil discovered that each local term $W_{\nu}(K, \chi)$ is best seen as the push-forward of a distribution $W_{\nu}(K_{\nu}, {\chi}_{\nu})$ from the multiplicative group of the $\nu -$adic field $K_\nu$ to $(0, \infty )$ (under $t_\nu \mapsto u = |t_\nu |_\nu)$. Indeed, all local terms, archimedean as non-archimedean, take then an (almost) identical functional form, and $W(K, \chi )$ is thus best seen as a distribution on the (classes of) ideles of $K$.

Some minor differences between our conventions and those of [We52]: our function $g$ is on the positive multiplicative group of the reals, not on the additive reals, and we do not shift by ${1 \over 2}$ the Mellin transform; we distribute the discriminant of the number field among its local components instead of keeping it as a multiple of the Dirac at $1$; we keep the poles and the zeros together on the same side of the Explicit Formulae.

Weil's ``Finite Part'' symbol hid away some differences between the finite and infinite places and in the case of a non-trivial character $\chi$, an important integral expression for the conductor exponents was needed, whose proof, perhaps because of its elementary nature, was in fact not spelled out (these integrals have an important extension to the non-abelian case [We74]).

The situation was improved by Haran [Ha90] who, for the Riemann Zeta function, gave to all Weil local terms an {\it exactly identical} formulation. He discovered that although they have at birth the shape of a {\it multiplicative} convolution, they also exist as an {\it additive} convolution. He chose to formulate his result in terms of Riesz Potentials, emphasizing the associated semi-group, but one can also just state it as\hfill\break\centerline{``The Explicit Formula is additive convolution with the Fourier Transform of $-\log(|y|)$''.}

In this paper, I show how to extend Haran's result to the general Dirichlet-Hecke $L$-series. The point of this is not so much the straightforward check that the Fourier Transform of $-\log(|y|)$ does the job in general, also at places of ramification, but rather to illustrate the additive-multiplicative tension mentioned above in a manner that treats all primes alike. It is also to be noted that Weil's local terms for archimedean places were not finite expressions, but involved limits and that finite expressions were apparently given only later ([Ba81]; of course this is related to the fact that the class of functions $g$'s considered in [We52] is much larger than the one considered here, which would allow for a very quick computation involving the partial fraction expansion of Euler's Gamma function, as for example in [Ha90]).

This derivation of the Explicit Formula gives birth to the ``conductor operator'':
$$H = \log(|x|_\nu) + \log(|y|_\nu)$$
I have used the physicist's notation so that $x$ and $y$ are Fourier dual variables. This operator acts as a hermitian (unbounded, but bounded below) operator on the Hilbert space $L^2(K_\nu, dx)$ of square-integrable functions on the local field $K_\nu$ with respect to the additive measure $dx$. It will be shown to be essentially self-adjoint on the domain consisting of the Bruhat-Schwartz functions, and that the local term of the Explicit Formula at the place $\nu$ actually is the {\it spectral analysis} of the conductor operator on $K_\nu$. This only works when the Explicit Formula is expressed as an integral over the critical line $Re(s) = {1\over2}$, which thus appears in a natural manner here.

The key feature allowing the spectral analysis of $H$ is that it is completely dilation invariant. Let $A$ be the operator of multiplication with $\log(|x|)$ and $B$ its Fourier conjugate which acts in the $y$-representation by multiplication with $\log(|y|)$. Let $U(t): \varphi(x) \mapsto {|t|}^{-1/2}\,\varphi({x\over t})$ be the unitary operators of dilation. Then $U(t)\,A\,U(t)^{-1} = A - \log(|t|)$, $U(t)\,B\,U(t)^{-1} = B + \log(|t|)$ so that $H$ is dilation invariant. And the same also holds for the commutator operator $K = i\,[B,\,A]$.

The complete dilation invariance of an operator $O$ means that it first needs to be transported to the multiplicative group $K_\nu^\times$ for its analysis, and should have there as generalized eigenvectors the multiplicative unitary characters $\chi(t)$, which from the additive side take the shape $\chi(x)\,|x|^{-1/2}$. Thus one expects an equation
$$O({\chi(x)}^{-1}|x|^{-s}) = O(\chi, s)\ {\chi(x)}^{-1}|x|^{-s}$$
for some ``spectral'' functions $O(\chi, s)$ on the critical line $Re(s) = 1/2$. One way to understand these equations is as an identity of distributions on $K_\nu$ (for the commutator operator this will be justified only on the {\it punctured} $\nu$-adics).

A most basic dilation invariant unitary is the composition of the additive Fourier Transform with the Inversion $\phi(x) \mapsto {1\over|x|}\phi({1\over x})$. Its spectral function is the Tate-Gel'fand-Graev Gamma function:
$$({\cal F}I)({\chi(x)}^{-1}|x|^{-s}) = \Gamma(\chi, s)\ {\chi(x)}^{-1}|x|^{-s}$$
Tate's Local Functional Equation states that this holds in the full critical strip $0 < Re(s) < 1$ as an identity of distributions on $K_\nu$.

The spectral function for the conductor operator is found to be:
$$H(\chi, s) = {\partial\over\partial s}\log(\Gamma(\chi,s))$$
and the identity makes sense for $s$ in the full critical strip but acquires its Hilbert space spectral interpretation only on the critical line. A similar result holds for the commutator operator $K = i\,[B, A]$:
$$K(\chi, s) = - i {\partial^2\over\partial s^2}\log(\Gamma(\chi,s))$$
which implies that $K$ is bounded.

It is found from this that the Inversion commutes with the conductor operator and anti-commutes with the commutator operator. Additional comments will be found in the ``Notes'' section concluding this paper.

{\bf B. Necessary facts from Tate's Thesis}

Let $K$ be a number field, $\AA$ its adele ring, $\AA^\times$ its multiplicative idele group, ${\cal C} = {\AA^\times}/{K^\times}$ its multiplicative group of idele classes.

On the local completion $K_\nu$ at the place $\nu$, one makes the following choice of basic additive character:
$$\lambda(x) = \exp(2\pi i\{\hbox{Sp}(x)\})$$
where Sp is the trace down to $\QQ_p$ (or minus the trace down to $\RR$ for an archimedean place), and $\{a\}$ the ``polar'' part of the $p-$adic number $a$. The Fourier Transform ${\cal F}(\varphi)$ of $\varphi$ (also denoted $\widetilde{\varphi}$) is defined through ${\cal F}(\varphi)(y) = \int_{K_\nu}\varphi(x)\lambda(-yx)\, dx$ where $dx$ is the unique self-dual additive Haar measure corresponding to the choice of $\lambda$ (so that ${\cal F}{\cal F}(\varphi)(x) = \varphi(-x)$). The $\nu$-adic module of $t \in K_\nu^\times$ is such that $d(tx) = |t|_\nu dx$. If the place is finite the module of a uniformizer is ${1\over q}$ with $q$ the cardinality of the residue field. The multiplicative Haar measure $d^*t$ on $K_\nu^\times$ is that constant multiple of $dt \over |t|_\nu$ such that, at a finite place, the subgroup of units has a unit volume, and, at an archimedean place, the push-forward to $(0, \infty)$ under the module is ${du\over u}$.

Let $\chi_\nu :K_\nu^\times \rightarrow U(1)$ be a local multiplicative unitary character. Let $\omega_s$ be the local principal character $t\mapsto |t|_\nu^s$, which is unitary for $s$ purely imaginary. For $Re (s) > 0$, $\chi_\nu(x)\omega_s(x)$ is locally integrable against ${dx\over|x|}$ hence defines a distribution on $K_\nu$ : 
$$\varphi_\nu \mapsto L_\nu (\chi_\nu\omega_s)(\varphi_\nu) = \int_{K_\nu^\times} \varphi_\nu(t\cdot 1) \chi_\nu(t)|t|_\nu^s\, d^*t$$
Up to a constant multiple this is the same as the distribution $\chi_\nu(x)|x|^{s-1}_\nu$. It has ([Ta50], see also [We66], [GeGrPi69]) a meromorphic continuation to all $s$, and satisfies:
$${\cal F}(\chi_\nu(x)|x|^{s-1}_\nu) = \Gamma_\nu(\chi_\nu, s)\chi_\nu^{-1}(y)|y|^{-s}_\nu$$
where the proportionality factor $\Gamma_\nu(\chi_\nu, s)$ is usually called the $\nu -$adic Tate-Gel'fand-Graev Gamma function. Tate tabulated these functions, and the notation separating $\chi_\nu$ and $s$ was introduced in [GeGrPi69]. As we will see this will be rather convenient for us. They are analytic and non-vanishing in the strip $0 < Re(s) < 1$, and periodic with period $2\pi i\over \log(q)$ if $\nu$ is finite. Tate also proved the following properties
$$\Gamma_\nu(\chi_\nu, s)\cdot\Gamma_\nu(\chi_\nu^{-1}, 1-s) = \chi_\nu(-1)$$
$$\overline{\Gamma_\nu(\chi_\nu, s)} = \chi_\nu(-1)\Gamma_\nu(\overline{\chi_\nu}, \overline{s})$$
$$Re(s) = {1\over2}\Rightarrow\left|\Gamma_\nu(\chi_\nu,s)\right|^2 = 1$$
The related function $\Lambda_\nu(\chi_\nu,s) = - {\Gamma^\prime_\nu(\chi_\nu,s)\over\Gamma_\nu(\chi_\nu,s)}$ satisfies
$$\Lambda_\nu(\chi_\nu,s) = \Lambda_\nu(\chi_\nu^{-1}, 1-s)$$
$$\overline{\Lambda_\nu(\chi_\nu,s)} = \Lambda_\nu(\overline{\chi_\nu},\overline{s})$$
$$Re(s) = {1\over2}\Rightarrow \Lambda_\nu(\chi_\nu,s) = \Lambda_\nu(\overline{\chi_\nu},\overline{s})\hbox{ and is real}$$

Let $F$ be a test-function on $\cal C$, i{.}e{.} a finite linear combination of functions $g(|t|){\chi(t)}^{-1}$ where $|t|$ is the global module, $g$ is a smooth function on $(0, \infty)$ with compact support, and $\chi$ is a unitary character on $\cal C$. The Mellin transform of such an $F$ is a function of quasi-characters (continuous homomorphisms $c: {\cal C} \rightarrow \CC^\times$) defined as
$$\widehat{F}(c)=\int_{\cal C}F(t)c(t)\, d^*_{\cal C}t$$
The Haar measure $d^*_{\cal C}t$ on $\cal C$ is normalized so that under push-forward under the module to $(0, \infty)$ it is sent to $du\over u$ (the fibers are compact).

All the Hecke $L$-functions $L(\chi, s)$ can be combined together as one meromorphic function $L(c)$ on the space of quasi-characters $c$, according to the formula $L(\chi\omega_s) = L(\chi, s)$ ($\omega_s(t)=|t|^s$). The multiplicity of $L$ at the quasi-character $c$ will be denoted by $m(c)$ (poles counted with positive multiplicities and zeros with negative multiplicities). Let:
$$Z(F)=\sum_{L(c)\ =\ 0\ or\ \infty}m(c)\widehat{F}(c)$$

Assuming from now on that $F(t) = g(|t|)\chi(t)^{-1}$, for a given fixed unitary character $\chi$, we see that $\widehat{F}(c)=\widehat{g}(s)$ for $c = \chi\omega_s$ while $\widehat{F}(c) = 0$ for all other quasi-characters $c$. Writing $Z(\chi, g)$ instead of $Z(F)$ we thus get its value as the sum of $\widehat{g}(s)$ over the zeros and poles of the $L$-function $L(\chi, s)$, as was considered in the Introduction.

It is convenient to consider the $L$-function $L(\chi\omega_s)$ not as valued in the complex numbers, but as valued in the distributions on the adeles \AA, in the following manner:
$$L(\chi\omega_s)(\varphi) = \int_{\AA^\times}\varphi(t\cdot 1)\chi(t)|t|^s\, d^*t$$
Here the test-function $\varphi(x)$ on the adeles is a finite linear combination of infinite products $\prod_\nu \varphi_\nu(x_\nu)$, where $\varphi_\nu$ depends only on the $\nu -$adic component, is for almost all $\nu$'s the characteristic function of the sub-ring of integers in $K_\nu$, is for all non-archimedean places locally constant with compact support (Bruhat function), and is for archimedean places a Schwartz function. Note that although this defines a distribution on \AA\ the integral is over $\AA^\times$, not over \AA. The Haar measure $d^*t$ (which lives on $\AA^\times$) is the (restricted) product of the local multiplicative Haar measures and it is not the same thing as $d^*_{\cal C}t$ above (which lives on $\cal C$).

In fact this defines directly $L(\chi\omega_s)(\varphi)$ only for $Re (s) > 1$, as an absolutely convergent product. Tate obtained the analytic continuation to all $s$ and proved the Functional Equation:
$$L(\chi\omega_s)(\widetilde{\varphi}) = L(\chi^{-1}\omega_{1-s})(\varphi)$$
Here $\widetilde{\varphi}(x)$ is the Fourier transform of $\varphi(x)$ so that the Fourier transform of the distribution $L(\chi\omega_s)$ is the distribution $L(\chi^{-1}\omega_{1-s})$.

{\bf C. The abelian Explicit Formula}

With a suitable test-function $\varphi = \prod_\nu \varphi_\nu$ the local factors of the Hecke $L$-function are obtained as
$$L_\nu(\chi_\nu, s) = L_\nu(\chi_\nu\omega_s)(\widetilde{\varphi_\nu})$$ (they have neither poles nor zeros in $Re(s) > 0$) while globally
$$L(\chi, s) = L(\chi\omega_s)(\widetilde\varphi)$$
and is thus obtained for $Re(s) > 1$ as the absolutely convergent product $\prod_\nu L_\nu(\chi_\nu, s)$. For $Re(s) < 0$ one represents $L(\chi, s)$ as the absolutely convergent product $\prod_\nu M_\nu(\chi^{-1}_\nu, 1-s)$, where  $M_\nu(\chi^{-1}_\nu, 1-s) = L_\nu(\chi^{-1}_\nu\omega_{1-s})(\varphi_\nu)$ (it has neither pole nor zero in $Re(s) < 1$). The Gamma function is
$$\Gamma_\nu(\chi_\nu, s) = {L_\nu(\chi_\nu, s) \over M_\nu(\chi^{-1}_\nu, 1-s)}$$

The computation of $Z(\chi, g)$ proceeds in the usual way of the calculus of residues as: $$Z(\chi, g) = {1 \over 2\pi i}\int_R \widehat{g}(s)\ -d\log\ L(\chi,s)$$ where R is the rectangle with corners $2-i\infty, 2 + i\infty, -1+i\infty, -1-i\infty$. For the necessary discussion allowing to discard the horizontal contributions, see [We52].

On the vertical line $Re(s) = 2$ one uses the product $\prod_\nu L_\nu (\chi_\nu, s)$ to express the integral as a sum of local terms. As the local $L$-functions have no poles or zeros in the half-plane $Re(s) > 0$, one can shift each local term to the line $Re(s) = c$, where $0 < c < 1$. (That $g$ has compact support turns out to imply that there are only finitely many non-vanishing terms, but this is not necessary to this argument). On the vertical line $Re(s) = -1$ one uses the product $\prod_\nu M_\nu (\chi_\nu^{-1}, 1-s)$ to express the integral as a sum of local terms. As the local integrands have no poles or zeros in the half-plane $Re(s) < 1$, one can shift each local term to the same line $Re(s) = c$.

Recombining the integrands from the left and the right, we end up with:
$$Z(\chi, g) = \sum_\nu W_\nu(\chi, g)$$
$$W_\nu(\chi, g) = {1 \over {2\pi i}}\int_{Re(s)=c} \widehat{g}(s)\Lambda_\nu(\chi, s)\, ds$$
$$\Lambda_\nu(\chi, s) = - {\partial\over\partial s}\log(\Gamma_\nu(\chi_\nu,s))$$

This differs from Weil's method in the way it uses the functional equation but it has the advantage of treating all places, ramified or not, alike. The next classical step from then on would be to obtain the inverse Mellin transform $w_\nu(\chi, u)$ of $\Lambda_\nu(\chi, s)$ (in the critical strip), which is a distribution on the positive half-line such that
$$W_\nu(\chi, g) = \int_0 ^\infty g({1 \over u})w_\nu(\chi, u)\, {du \over u}$$
and then (somewhat surprisingly in the real case and quite surprisingly in the complex case) to realize that it comes from $K_\nu^\times$ through $t \mapsto u = |t|_\nu$, and to guess from the functional form of the final archimedean results the way to deal with the conductor of the character.

Let us now work locally and drop all subscripts $\nu$ to lighten the notation. The identities to follow are identities of locally integrable functions, taken in the distributional sense, and the parameter $s$ is in the critical strip $0<Re(s)<1$.
$${\cal F}(\chi(y)|y|^{s-1}) = \Gamma(\chi, s)\chi^{-1}(x)|x|^{-s}$$
$${\cal F}\bigl(\log|y|\,\chi(y)|y|^{s-1}\bigr)=\Gamma^\prime(\chi,s)\chi^{-1}(x)|x|^{-s} + \Gamma(\chi, s)\,(-\log|x|)\chi^{-1}(x)|x|^{-s}$$
$$\chi^{-1}(x)\Lambda(\chi,s)|x|^{-s}=-\log|x|\,\chi^{-1}(x)|x|^{-s} + {\cal F}\left(-\log|y|\cdot\chi(y){|y|^{s-1}\over\Gamma(\chi, s)}\right)$$
We think of these identities as being applied to some given test-function on the $\nu-$adic field and we now want to apply a further operation of integration, this time a complex line integral against ${1 \over {2\pi i}}\int_{Re(s)=c} \widehat{g}(s)\ \, ds$, for $0<c<1$. One needs a lemma according to which (pointwise in $y \neq 0$)
$${1 \over {2\pi i}}\int_{Re(s)=c} \widehat{g}(s)\chi(y){|y|^{s-1}\over\Gamma(\chi, s)}\, ds
= {\cal F}^{-1}(\chi^{-1}(x)g(|x|))(y)$$
The integral converges absolutely due to the rapid decrease of $\widehat{g}(s)$ whereas $\left|{1\over\Gamma(\chi,s)}\right|= \left|\Gamma(\chi^{-1},1-s)\right|$ ($=1$ on the critical line) is in any case periodic for a finite place, and controlled by the Stirling formula to be $O(|s|^{1/2-c})$ for a real place, $O(|s|^{1-2c})$ for a complex place. So the left-hand-side defines a continuous function of $y$ (for $y\neq0$) and also a tempered distribution as it is $O(|y|^{c-1})$. Applying both sides to a test function we see  (using Mellin inversion) that they coincide as distributions hence pointwise (for $y\neq0$).

As everything converges absolutely we are allowed to intervert integrals. Defining $$W_\nu(\chi, g; x) = {1 \over {2\pi i}}\int_{Re(s)=c} \widehat{g}(s)\Lambda_\nu(\chi, s)|x|_\nu^{-s}\chi_\nu^{-1}(x)\, ds$$
and $F_\nu(x) = g(|x|_\nu)\chi_\nu^{-1}(x)$ ($F_\nu(0) = 0$) we obtain (for $x \neq 0$):
$$W_\nu(\chi, g; x) = -\log(|x|_\nu) F_\nu(x) + {\cal F}\left(-\log(|y|_\nu)\cdot{\cal F}^{-1}(F_\nu)(y)\right)$$
$$W_\nu(\chi, g; x) = -\log(|x|_\nu) F_\nu(x) + {\cal F}\left(-\log(|y|_\nu)\right) * F_\nu(x)$$
where $*$ is the symbol of additive convolution on $K_\nu$ . Plugging in $x = 1$ gives:
$$W_\nu(F) = W_\nu(\chi, g) = (G_\nu * F_\nu)(1)$$
where $G_\nu$ is the Fourier Transform of $-\log(|y|_\nu)$ and $*$ is an additive convolution.

{\bf Theorem C1}: For an arbitrary test-function $F$ on the idele classes (a finite linear combination of smooth compactly supported functions of the global module twisted with unitary idele class characters), with $F_\nu$ the restriction of $F$ to $K_\nu$ (taking the value $0$ at the origin), and with $G_\nu(x) = {\cal F}(-\log(|y|_\nu))$
$$\sum_{L(c)\ =\ 0\ or\ \infty}m(c)\widehat{F}(c) = \sum_\nu(G_\nu * F_\nu)(1)$$

{\bf Note}: for the rational number field and trivial character, this is due to Haran [Ha90], in a distinct but directly equivalent formulation.

{\bf D. Fourier transform of $-\log (|y|)$}

To put the result in Weil's form we need to obtain the Fourier transform of $-\log(|y|)$. This is a classical calculus exercise in the archimedean case (homogeneity alone immediately gives the result away from the origin, but we want the exact determination). Over $\QQ_p$ it has been obtained by Vladimirov [Vl88]. It is also implicitely contained in Haran's work. As previously for the Explicit Formula we derive it in a simple manner from the properties of homogeneous distributions considered in [We66], [GeGrPi69].

Again we work locally, and will drop most $\nu$ subscripts to ease up on the notations. Let $\Delta_s$ be the homogeneous distribution, given for $Re(s) > 0$, by $\Delta_s(\varphi) = \int_{K_\nu}\varphi(x)|x|^{s-1}\, dx$. This time we need an explicit expression for its analytic continuation around $s = 0$. Let us (temporarily) pick a specific function $\omega(x)$ as follows:

If $\nu$ is finite, we take $\omega$ to be the characteristic function of the ring of integers $\cal O$. Let $\pi$ be a uniformizer and $q = {|\pi|}^{-1}$ (which is also the cardinality of the residue field). Let $\delta$ be the differental exponent at $\nu$, which is the largest integer such that $(x \in \pi^{-\delta}{\cal O}) \Rightarrow (\hbox{Sp}(x) \in \ZZ_p)$. Then Vol$({\cal O}) = q^{-\delta/2}$, and:
$$\Delta_s(\omega) = q^{-\delta/2} {1 - {1\over q} \over 1 - {1\over q^s}}$$
So $\Delta_s(\omega)$ has a meromorphic continuation to all $s$, and around $0$:
$$\Delta_s(\omega) = {R\over s} + P(\omega) + \cdots$$
with $R = q^{-\delta/2} {1 - {1\over q} \over \log(q)}$ and $P(\omega) = {1\over 2}\log(q)\cdot R$

If $\nu$ is real, we take $\omega(x) = \exp(- \pi x^2)$ and then (here $\Gamma$ is Euler's Gamma function):
$$\Delta_s(\omega) = \pi^{-s/2}\, \Gamma(s/2) = {R\over s} + P(\omega) + \cdots$$
with $R = 2$ and $P(\omega) = -\log(\pi) - \gamma_e$ ($\gamma_e = -\Gamma^\prime(1)$ = Euler-Mascheroni's constant)

If $\nu$ is complex we take $\omega(z) = \exp(- \pi z\overline z)$ and recall that $|z|$ means $z\overline z$, and that the self-dual Haar measure for $\lambda(z) = \exp(-2\pi i\ 2Re(z))$ is twice the Lebesgue measure, so that:
$$\Delta_s(\omega) = 2\pi\ \pi^{-s}\, \Gamma(s) = {R\over s} + P(\omega) + \cdots$$
with $R = 2\pi$ and $P(\omega) = - 2\pi(\log(\pi) + \gamma_e)$

The continuation to $Re(s) > -1$ (to all $s$ if $\nu$ is non-archimedean) of $\Delta_s$ is obtained (for example) by:
$$\Delta_s(\varphi) = \Delta_s(\omega)\varphi(0) + \int_{K_\nu}(\varphi(x)-\varphi(0)\omega(x))\, |x|^s\, {dx \over |x|}$$
Hence $\Delta_s$ has a simple pole at $0$ ($\delta_0$ = Dirac at $0$, not to be confused with the differental exponent!):
$$\Delta_s = {{R}\over s}\cdot \delta_0 + P + \cdots$$
$$P(\varphi) = P(\omega)\varphi(0) + \int_{K_\nu}(\varphi(x)-\varphi(0)\omega(x))\, {dx \over |x|}$$
Adopting the notation $P_\omega$ for the distribution $\int_{K_\nu}(\varphi(x)-\varphi(0)\omega(x))\, {dx \over |x|}$, now for an arbitrary test-function $\omega$ such that $\omega(0) = 1$ we thus see that the ``intrinsic finite part'' $P$ can be expressed as $$P(\varphi) = P_\omega(\varphi) + P(\omega)\delta_0(\varphi)$$

We know that ${\cal F}(\Delta_s(y)) = \Gamma(s)\Delta_{1-s}(x)$ with $\Gamma(s)$ the Gamma function for the trivial character. Expanding this in powers of $\varepsilon$ with $s = 1 - \varepsilon$, we get, on the left hand side:
$${\cal F}(\Delta_{1-\varepsilon}) =  \delta_0 + \varepsilon{\cal F}(-\log|y|) + \cdots$$
On the right hand side, we have, introducing constants $\gamma$ and $\tau$:
$$\Gamma(1 - \varepsilon) = \gamma(\varepsilon + \tau\varepsilon^2 + \cdots)$$
$$\Delta_\varepsilon = {R\over\varepsilon}\cdot\delta_0 + P + \cdots$$
Hence $\gamma R = 1$ and the formula for $G(x) = {\cal F}(-\log|y|)$ is:

{\bf Theorem D1}\quad$G = \gamma P + \tau\delta_0 = \gamma P_\omega + (\gamma P(\omega) + \tau)\delta_0 = \gamma P_\omega + G(\omega)\delta_0$

$G(\varphi) = \gamma P(\varphi) + \tau \varphi(0)$ is also directly given by:
$$G(\varphi) =  \int_{K_\nu}(-\log|y|)\, \widetilde{\varphi}(y)\, dy = - \left.{\partial \over {\partial s}}\right|_{s = 1}\ \Delta_s(\widetilde{\varphi})$$
from which follows $G(\varphi_t) = G(\varphi) - \log(|t|)\varphi(0)$ for $\varphi_t(x) = \varphi(tx)$.

For all places, the constant $\gamma = 1/R = -\Gamma^\prime(1)$ turns out to be such that the multiplicative Haar measure $d^\times t = \gamma {dt \over |t|}$ on $K_\nu^\times$ has the property  $\int_{1\leq|x|<|X|}\, d^\times t = \log(|X|)$ (with $X \in K_\nu, |X| \geq 1)$. For a discrete $\nu$ it assigns a value of $\log(q)$ to the volume of the units ($d^\times t = \log(q) d^*t$).

Let us suppose now that the place $\nu$ is non-archimedean. For $\omega$ the characteristic function of the integers, its Fourier transform is $\widetilde\omega(y) = q^{-\delta/2}\omega(\pi^\delta y)$ and $\Delta_s(\widetilde\omega) = q^{\delta(s-1)}{1 - {1\over q} \over 1 - {1\over q^s}}$ so that $G(\omega) = -\left.{\partial \over {\partial s}}\right|_{s = 1}\ \Delta_s(\widetilde{\omega}) = {\log(q)\over q-1} - \delta\log(q)$. So:
$$G(\varphi)= \gamma\cdot P_\omega(\varphi) + \gamma\cdot q^{-{\delta\over2}-1}\varphi(0) - \delta\log(q)\varphi(0)$$

Let us apply this to the characteristic function $\psi$ of the ball of radius $1/q$ centered at $1$: we get $G(\psi) = \gamma\cdot q^{-{\delta\over2}-1}$. So it is advantageous to replace $\omega$ with $\omega_1 = \omega - \psi$, and then:
$$G(\varphi) = \gamma P_{\omega_1}(\varphi) - \delta\log(q)\varphi(0)$$
Going back to Weil's local term $W_\nu(F) = (G_\nu * F_\nu)(1)$, this gives:
$$W_\nu(F) =  \int_{K_\nu}(F_\nu(1-x) - F_\nu(1)\omega_1(x))\, \gamma{dx \over |x|} - \delta\log(q) F_\nu(1)$$
and after the change of variable $x \mapsto 1 - 1/t$, we get:

{\bf Theorem D2 ([We52])}
$$W_\nu(F) = \int_{|t| = 1}(F_\nu({1\over t}) - F_\nu(1))\, {d^\times t \over |1 - t|} + \int_{|t| \neq 1}F_\nu({1\over t})\, {d^\times t \over |1 - t|}-\ \delta\log(q) F_\nu(1)$$

This is Weil's local term, with the addition of the local component of the discriminant. The integral over $|t| = 1$ for the test-function $F(t) = f(|t|)\chi^{-1}(t)$ gives a non-zero result only if $\chi$ is ramified at $\nu$, whereas the other integral gives a non-zero result only when $\chi$ is non-ramified at $\nu$. We will come back to the ramified characters, but first let us examine the local term at an archimedean spot.

First we consider a real spot $\nu$.\hfil\break For $\omega(x) = \exp(- \pi\ x^2)$ we evaluated $\Delta_s(\omega) = \pi^{-s/2}\, \Gamma(s/2)$, $R = 2$ (so $\gamma = 1/2$), and $P(\omega) = -(\log(\pi) + \gamma_e)$. Using $\omega = \widetilde\omega$  (and also $\Gamma^\prime(1/2) = -\ \Gamma(1/2)(\gamma_e + 2\ \log(2)) )$ we get $G(\omega) = - \left.{\partial \over {\partial s}}\right|_{s = 1}\ \Delta_s(\widetilde{\omega}) = (\log(\pi) + \gamma_e + 2\log(2))/2$. From this we get the constant $\tau = G(\omega) - \gamma P(\omega) = \log(2\pi) + \gamma_e$.

It is convenient now to choose other functions $\omega$ (always obeying the condition $\omega(0) = 1$). For example, we can take $\omega(x) = 1$ if $|x| < 1$, $\omega(x) = 0$ if $|x| \geq 1$. (it doesn't matter that $\omega$ is not smooth away from $0$). For this $\omega$, $\Delta_s(\omega) = 2/s$, so $P(\omega) = 0$ and $G(\omega) = \tau = \log(2\pi) + \gamma_e$. So we get:
$${\cal F}(-\log|y|)(\varphi) = \int_{|x| \leq 1}(\varphi(x) - \varphi(0))\, {dx \over 2|x|} + \int_{|x| > 1}\varphi(x)\, {dx \over 2|x|}+\ (\log(2\pi) + \gamma_e)\varphi(0)$$
This gives the following expression for the Weil term ($d^\times t = {dt \over 2|t|}$):

{\bf Theorem D3}
$$W_\RR(F) = \int_{1/2}^\infty(F({1\over t}) - F(1))\, {d^\times t \over |1 - t|} + \int_{-\infty}^{1/2}F({1\over t})\, {d^\times t \over |1 - t|} +\ (\log(2\pi) + \gamma_e) F(1)$$

Choosing $\omega(x) = \sin(\pi x)/\pi x$ whose Fourier transform (as a distribution) is the characteristic function of $-1/2 \leq x \leq 1/2$, we get $\Delta_s(\widetilde\omega) = {1\over s}{1\over 2^{s-1}}$ hence $G(\omega) = 1 + \log(2)$ and:
$${\cal F}(-\log|y|)(\varphi) =  \int_\RR \left({\varphi(x) - \varphi(0){\sin(\pi x) \over \pi x}}\right)\, {dx \over 2|x|} + (1 + \log(2))\varphi(0)$$
This leads to the following expression for Weil's term:
$$W_\RR(F) = \int_{\RR^\times}\left({F({1\over t}) - {{{t\over\pi}\, \sin({\pi\over t})}\over t - 1}\, F(1)}\right)\, {d^\times t\over|1 - t|} + (1 + \log(2)) F(1)$$
As $\exp(G(\omega)) = 2e$, the function $\omega(x) = \sin(2\pi ex)/2\pi ex$ would lead to a similar expression with no Dirac term.

Choosing $\omega(x) = \left({\sin(\pi x)/\pi x}\right)^2$ whose Fourier transform is the triangle function $x \mapsto 1 - |x|$ for $|x| \leq 1$, we get $\Delta_s(\widetilde\omega) = {2\over s(s+1)}$ hence $G(\omega) = 3/2$ and:
$${\cal F}(-\log|y|)(\varphi) =  \int_\RR \left({\varphi(x) - \varphi(0){\left({\sin(\pi x) \over \pi x}\right)^2}}\right)\, {dx \over 2|x|} + (3/2)\varphi(0)$$
$$W_\RR(F) = \int_{\RR^\times}\left({F({1\over t}) - {\left({{{t\over\pi}\, \sin({\pi\over t})}\over t - 1}\right)^2}\, F(1)}\right)\, {d^\times t\over|1 - t|} + (3/2) F(1)$$

Switching to a complex spot $\nu$, we obtained with $\omega(z) = \exp(-\ \pi z\overline z)$: $\Delta_s(\omega) = 2\pi\ \pi^{-s}\, \Gamma(s)$, $R = 2\pi$, $\gamma = 1/2\pi$ and $P(\omega) = -2\pi(\log(\pi) + \gamma_e)$. The Fourier transform (recall that $dz = 2rdrd\theta$ and $|z| = z\overline z$) of $\omega(z)$ is $2\omega(2z)$, so $\Delta_s(\widetilde \omega) = (4\pi)^{1-s}\ \Gamma(s)$ and we get $G(\omega) = - \left.{\partial \over {\partial s}}\right|_{s = 1}\ \Delta_s(\widetilde{\omega}) =  \log(4\pi) + \gamma_e$. From this we get the constant $\tau = G(\omega) - \gamma P(\omega) = 2(\log(2\pi) + \gamma_e)$. Choosing $\omega(z)$ to be the characteristic function of the unit disc, we find $\Delta_s(\omega) = 2\pi/s$ hence $P(\omega) = 0$ and so $G(\omega) = \tau$, and the following thus emerges:
$${\cal F}(-\log|y|)(\varphi) = \int_{|x| \leq 1}(\varphi(x) - \varphi(0))\, {dx \over 2\pi|x|} + \int_{|x| > 1}\varphi(x)\, {dx \over 2\pi|x|}+\ 2(\log(2\pi) + \gamma_e)\varphi(0)$$
This gives the following expression for the Weil term ($d^\times t = drd\theta/\pi r$):

{\bf Theorem D4}
$$W_\CC(F) = \int_{Re(t)\geq 1/2}(F({1\over t}) - F(1))\, {d^\times t \over |1 - t|} + \int_{Re(t)\leq 1/2}F({1\over t})\, {d^\times t \over |1 - t|}+\ 2(\log(2\pi) + \gamma_e) F(1)$$

There are of course countless other possible choices for the function $\omega$ and ensuing expressions for $W_\RR(F)$ and $W_\CC(F)$.

{\bf E. The Explicit Formula and the conductor operator\par}

In the remaining sections of this paper we will only consider one place $\nu$ of the number field $K$ at a time, and we will drop most $\nu$-subscripts from our notations. Let us first return to the local contribution to the Explicit Formula at finite place. Let $\chi$ be a local, ramified, unitary character (ramified meaning that the restriction of $\chi$ to ${\ |t| = 1\ }$ is non-trivial). Let $f\geq1$ be the conductor exponent of $\chi$, that is the smallest positive integer such that $| 1 - t | \leq q^{-f} \Rightarrow \chi(t) = 1$. For $q=2$ the smallest possible value of $f$ is $2$ whereas for all other $q$'s it is $1$. The Gamma function of $\chi$ was computed by Tate to be simply $q^{(f+\delta) s}$ up to a certain multiplicative constant, not important here. So minus its logarithmic derivative is the constant $-(f + \delta) \log(q)$.

Going back to the evaluation of $$W_\nu(\chi, g; x) = {1 \over {2\pi i}}\int_{Re(s)=c} \widehat{g}(s)\Lambda_\nu(\chi, s)|x|_\nu^{-s}\chi_\nu^{-1}(x)\, ds$$
this gives $W_\nu(\chi, g; x) = -(f + \delta) \log(q) g(|x|) \chi^{-1}(x)$, and for $F(x) = g(|x|)\chi^{-1}(x)$ our previous computation of $W_\nu(\chi, g; x)$ then gives:
$$-(f + \delta) \log(q) F(x) = (-\log|x|) F(x) + {\cal F}(-\log|y|) * F(x)$$
Choosing the function $g$ to have its support contained in $1/q < u < q$, and such that $g(1) = 1$, plugging in $x = 1$, and using our evaluation of ${\cal F}(-\log|y|)$, we end up with:$$-(f + \delta) \log(q) =  \int_{|t| = 1}(\chi(t) - 1){{d^\times t}\over{|1-t|}} - \delta\log(q)$$

{\bf Theorem E1 ([We74])} $f \log(q) =\int_{|t| = 1}(1 - \chi(t)){{d^\times t}\over{|1-t|}}$

It is not difficult to evaluate directly this integral which is a necessary component of the final result of [We52] (perhaps for this reason a proof was omitted there).

We now turn to an operator theoretic interpretation of the general local contribution to the Explicit Formula, also at archimedean places. For this let $L^2$ be the Hilbert Space of square integrable functions on $K_\nu$ with respect to the additive Haar measure $dx$. Let ${\cal S}$ be the domain of Schwartz-Bruhat functions on $K_\nu$, and ${\cal S}_0$ the subdomain consisting of functions having compact support and vanishing identically in a neighborhood of the origin. We defined in the Introduction the operators $A: \varphi(x) \mapsto \log(|x|)\cdot\varphi(x)$, $B = {\cal F}A{\cal F}^{-1} = {\cal F}^{-1}A{\cal F}$. The conductor operator is defined on ${\cal S}$ as follows:
$$H(\varphi)(x) = (A+B)(\varphi)(x) = \log(|x|)\,\varphi(x) + {\cal F}\left(\log(|y|)\,{\cal F}^{-1}(\varphi)(y)\right)(x)$$
The right-hand side is at any rate a continuous function which belongs to $L^2$. There is also the Inversion $I: \varphi \mapsto (x \mapsto {1\over |x|}\varphi({1\over x}))$ which acts on ${\cal S}_0$ and (isometrically) on $L^2$.

{\bf Theorem E2:} The conductor operator $H$ with initial domain ${\cal S}_0$ is essentially self-adjoint. It is unbounded, but bounded below. It commutes with the unitary dilations, the Fourier transform and the Inversion.

Let also $K = i[B, A]$.

{\bf Theorem E3:} The commutator operator $K$ with initial domain ${\cal S}_0$ is essentially self-adjoint. It is bounded. It commutes with the unitary dilations and anti-commutes with the Fourier Transform and the Inversion.

The remaining sections are devoted to developments including the completion of the proofs of Theorems E2 and E3. From the point of view of the Hilbert Space, the crucial thing is of course the proof of the self-adjointness of $(H, {\cal S}_0)$ and the determination of the full domain ${\cal D}(H)$ of $H$. But this we obtain only at the very end. In fact the most interesting domain is not ${\cal S}$, ${\cal S}_0$ or ${\cal D}(H)$ but a fourth one $\Delta$ that has the property of being stable under ${\cal F}, A, B$, and $I$. Statements such that ``$H$ commutes with the Inversion'' are to be understood in the mean time to apply to test-functions in ${\cal S}_0$.

{\bf F. Concrete spectral analysis on $\QQ_p$}

Let us first elucidate the spectrum of the conductor and commutator operators over $\QQ_p$ using very concrete methods. As $H$ and $K$ commute with averaging over the units one can split the discussion between the ramified and the non-ramified spectrum. 

First the ramified spectrum. From the above discussion
$$H(\varphi)(x) = f\log(p)\,\varphi(x)$$
if $\varphi$ has the homogeneity of a ramified character with conductor exponent $f$. So $HA = AH$ on such $\varphi$'s and $BA = (H-A)A = A(H-A) = AB$ so that $K(\varphi) = 0$.

{\bf Theorem F1:} The ramified spectrum of the conductor operator is discrete. It has infinite multiplicity and support $\{\log(p), 2\log(p), 3\log(p), \dots\}$ if $p>2$, $\{2\log(2), 3\log(2), \dots\}$ if $p=2$. The commutator operator vanishes on the ramified subspace of $L^2$.

We can now restrict $H$ to the invariant part $L^2_{inv}(\QQ_p)$ spanned by functions depending only on the module $|x|$. Let's first evaluate $H(\theta)$ for $\theta$ the characteristic function of the units. For this one can use the Fourier transform $G(x)$ of $-\log(|y|)$ determined above:
$$G(\varphi) = {\log(p) \over 1 - 1/p}\left(\int_{|x|\leq 1}(\varphi(x) - \varphi(0))\, {dx\over|x|} + \int_{|x|>1}\varphi(x)\, {dx\over|x|} + {1\over p}\varphi(0)\right)$$
and the outcome for $H(\theta)(x) = \log(|x|)\theta(x) - (G*\theta)(x)$ is
$$\eqalign{
(|x| < 1)\quad\ H(\theta)(x) &= -\log(p) \cr
(|x| = 1)\quad\ H(\theta)(x) &= \;0 \cr
(|x| > 1)\quad\ H(\theta)(x) &= -{\log(p)\over|x|} }$$
By dilation invariance one has $H(\theta_i)(x) = H(\theta)(p^ix)$ for $\theta_i(x) = \theta(p^ix)$. So the matrix elements of $H$ in the orthogonal basis $\{\theta_i, i\in\ZZ\}$ depend only on $i-j$. An orthonormal basis of $L^2_{inv}(\QQ_p)$ is given by the functions
$$\eta_i(x) = (1 - 1/p)^{-1/2} p^{-i/2} \theta(p^i x)$$
and the matrix elements of $H$ in this basis are
$$\eqalign{
(i\neq j)\quad\ H_{ij} &=\ <H\eta_i|\eta_j>\ = -\log(p)\,p^{-|i-j|/2} \cr
(i = j)\quad\ H_{ii} &=\ <H\eta_i|\eta_i>\ = 0 }$$
We can now use the isometry from $L^2_{inv}(\QQ_p)$ to $L^2(S^1, {d\theta\over2\pi})$ which sends $\eta_i$ to $e_i(z) : z\mapsto z^i$. Under this isometry $H$ is transported to the multiplication operator with the function $$h(z) = -\log(p)\,\left({z\over\sqrt{p}-z}+{\overline{z}\over\sqrt{p}-\overline{z}}\right)$$

{\bf Theorem F2:} The conductor operator is bounded on $L^2_{inv}(\QQ_p)$ and has a continuous spectrum with support $[{-2\log(p)\over\sqrt{p}-1}, {+2\log(p)\over\sqrt{p}+1}]$.

Let's now evaluate $K(\theta) = i\,(BA - AB)(\theta)$. As $A(\theta) = 0$ this is just $-i\,A(H(\theta))$ so 
$$\eqalign{
(|x| < 1)\quad\ K(\theta)(x) &= i\,\log(p) \log(|x|)\cr
(|x| = 1)\quad\ K(\theta)(x) &= \;0 \cr
(|x| > 1)\quad\ K(\theta)(x) &= i\,\log(p){\log(|x|)\over|x|} }$$

The matrix elements of $K$ in the orthonormal basis $\{\eta_j, j\in\ZZ\}$ depend only on $j-k$ and are
$$\eqalign{
(j\neq k)\quad\ K_{jk} &=\int_{\QQ_p}\overline{K\eta_j(x)}\eta_k(x)\,dx\ = -i\,(\log(p))^2\,(k-j)\,p^{-|k-j|/2} \cr
(j = k)\quad\ K_{jj} &= 0 }$$
Hence $K$ is transported to the multiplication operator with the function $$k(z) = i\log(p)^2\ \sum_{k\neq0}k\,p^{-|k|/2}e_k(z)$$
This is $-i\,\log(p)\,z{\partial\over\partial z}h(z)$. So
$$k(z) = i\,\log(p)^2\,\left({{z\over\sqrt{p}}\over{(1 - {z\over\sqrt{p}})^2}}-{{\overline{z}\over\sqrt{p}}\over{(1 - {\overline{z}\over\sqrt{p}})^2}}\right)$$
A little calculation then completes the proof of the following statement:

{\bf Theorem F3:} The commutator operator is bounded on $L^2_{inv}(\QQ_p)$ and has a continuous spectrum with support the closed interval $[-{2\log(p)^2\over\sqrt{p}},+{2\log(p)^2\over\sqrt{p}}]$.

Before reworking the previous paragraphs into more intrinsic terms, I will make an apart\'e on one approach, certainly very well-known, to the Poisson summation formula.

Let $u >1$, $L_u(s) = 1/(1-u^{-s})$, $\Lambda_u(s) = - L^\prime_u(s)/L_u(s) = \log(u)\,u^{-s}/(1-u^{-s})$.
$$Re(s) > 0\;\Rightarrow\;\Lambda_u(s) = \log(u)\,\sum_{k\geq1}u^{-ks}$$
$$\Lambda_u(s) + \Lambda_u(-s) = -\log(u)$$
For a smooth function $g$ with compact support on $(0,\infty)$ the first of these equations and Mellin inversion implies
$${1\over{2\pi i}}\int_{Re(s)=1} \widehat{g}(s)\Lambda_u(s)\, ds = \log(u)\,\sum_{k>0}g(u^k)$$
while the second equation (functional equation) gives
$$\eqalign{
{1\over{2\pi i}}\int_{Re(s)=-1} \widehat{g}(s)\Lambda_u(s)\, ds &= -{1\over{2\pi i}}\int_{Re(s)=1} \widehat{g}(-s)(\Lambda_u(s)+\log(u))\, ds \cr
&= -\log(u)g(1) - \log(u)\,\sum_{k<0}g(u^k)
}$$
The poles of $\Lambda_u(s)$ are at $j\,{2\pi i\over\log(u)}, j\in\ZZ$ with residue $+1$ so that one obtains in the end
$$\log(u)\,\sum_{k\in\ZZ}g(u^k) = \sum_{j\in\ZZ}\widehat{g}(j\,{2\pi i\over\log(u)})$$
which is one formulation of Poisson summation formula.

{\bf G. Transformation Theory and an intrinsic spectral analysis\par}

We go back to the general number field $K$ and its completion $K_\nu$ (archimedean or finite). We chose the multiplicative Haar measure $d^*t$ on $K_\nu^\times$ to assign a volume of $1$ to the units if $\nu$ is finite, and to be mapped to $du/u$ under $t\mapsto u = |t|$ for an archimedean place. For the positive constant $\alpha$ such that $d^*t = \alpha^{-2}\,{dt\over |t|}$ the maps
$$\varphi(x)\mapsto\{t\mapsto\alpha\sqrt{|t|}\varphi(t)\}$$
$$f(t)\mapsto\{x\mapsto\alpha^{-1}|x|^{-1/2}f(x)\}$$
establish a unitary isomorphism between $L^2(K_\nu,dx)$ and $L^2(K_\nu^\times,d^*t)$.

If we transport the conductor operator from $L^2(K_\nu,dx)$ to $L^2(K_\nu^\times,d^*t)$ we find that its dilation invariance becomes a commutation property with the multiplicative translations on $K_\nu^\times$. Hence, by the well-known generalization of Fourier Theory, $H$ can be completely analyzed in terms of the (unitary) characters of that locally compact abelian group. Although it is naturally defined in additive terms, its analysis proceeds in multiplicative terms. So let $X_\nu$ be the dual of $K_\nu^\times$, that is the locally compact abelian group of unitary multiplicative characters. For a finite place $X_\nu$ is the union of countably many circles indexed by the characters of the compact subgroup of units, for the real place it is the union of two real lines indexed by the trivial and the sign characters, for the complex place it is the union of countably many real lines indexed by the characters of the unit circle.

There is a unique Haar measure $\mu$ on $X_\nu$ so that the inversion formula of Fourier Theory and then the Plancherel identity hold:
$$\eqalign{
f(\chi) &= \int_{K_\nu^\times}f(x)\chi(t)\,d^*t \cr
f(t) &= \int_{X_\nu}f(\chi)\,\overline{\chi(t)}\,d\mu(\chi) \cr
\int_{K_\nu^\times}|f(t)|^2\,d^*t &= \int_{X_\nu}|f(\chi)|^2\,d\mu(\chi)
}$$

Under these isometries the Inversion $\varphi(x)\mapsto{1\over|x|}\varphi({1\over x})$ on $L^2(K_\nu,dx)$ ($f(t)\mapsto f(1/t)$ on $L^2(K_\nu^\times,d^*t)$) becomes the operation of complex conjugation $f(\chi) \mapsto f(\overline{\chi})$ (on the underlying space).
It is also interesting to look at the fate of the additive Fourier Transform ${\cal F}$. First transported to $L^2(K_\nu^\times,d^*t)$ it maps $f(t)$ (``smooth'' with compact support) to $F(t) = \sqrt{|t|}\,{\cal F}(|t|^{-1/2}f(t))$. This implies:
$$F(\chi) = \int_{K_\nu^\times} F(t)\chi(t)\,d^*t = \alpha^{-2}\ \int_{K_\nu} {\cal F}(|x|^{-1/2}f(x))\,\chi(x)|x|^{-1/2}\,dx $$
$$F(\chi) = \alpha^{-2}\ \Gamma(\chi,1/2)\ \int_{K_\nu} |x|^{-1/2}f(x)\,\chi(1/x)|x|^{-1/2}\,dx$$
$$F(\chi) = \Gamma(\chi,1/2)\ f(\overline{\chi})$$

We have obtained the first statements of:

{\bf Theorem G1:} Under the isometries $L^2(K_\nu,dx) \sim L^2(K_\nu^\times,d^*t) \sim L^2(X_\nu,d\mu)$ the inversion is transported to the unitary $$f(\chi)\mapsto f(\overline{\chi})$$ and the additive Fourier transform is transported to the unitary
$$f(\chi) \mapsto F(\chi) = \Gamma(\chi,1/2)\ f(\overline{\chi})$$ while the conductor operator is transported to the hermitian operator
$$f(\chi) \mapsto F(\chi) = - \Lambda(\chi,1/2)\ f(\chi)$$

{\bf Corollary:} The conductor operator commutes with the inversion.

{\bf Proof (of Corollary):} $\Lambda(\chi,1/2) = \Lambda(\overline{\chi},1/2)$.

{\bf Proof (of Theorem):}
In our derivation of the Explicit Formula we obtained the following: let $g$ be a smooth function with compact support on $(0,\infty)$ and Mellin Transform
$$\widehat{g}(s) = \int_{0}^{\infty} g(u){u}^s \, { du \over u}$$
and let $\chi$ be a unitary character on $K_\nu^\times$. Then with $\varphi(x) = g(|x|) \chi^{-1}(x)$ for $x\neq0$, $\varphi(0) = 0$, one has for $x\neq0$
$$H(\varphi)(x) = -{1 \over {2\pi i}}\int_{Re(s)=1/2} \widehat{g}(s)\Lambda(\chi, s)|x|^{-s}\chi^{-1}(x)\, ds$$

Let's write now $f(t) = \sqrt{|t|} g(|t|)\chi^{-1}(t)$ so that the conductor operator transported to $K_\nu^\times$ sends $f$ to
$$F(t) = -{1 \over {2\pi i}}\int_{Re(s)=1/2} \widehat{g}(s)\Lambda(\chi, s)|t|^{{1\over2}-s}\chi^{-1}(t)\, ds$$
$$F(t) = -{1 \over {2\pi}}\int_{-\infty}^{+\infty} \widehat{g}({1\over2} + i\tau)\Lambda(\chi, {1\over2} + i\tau)|t|^{-i\tau}\chi^{-1}(t)\, d\tau$$
Let's write $\omega_\tau(t) = |t|^{i\tau}$. Then
$$F(t) = -{1 \over {2\pi}}\int_{-\infty}^{+\infty} \widehat{g}({1\over2} + i\tau)\Lambda(\chi\omega_\tau, {1\over2})\,\overline{\chi\omega_\tau(t)}\, d\tau$$
We now split the discussion according to whether $\nu$ is finite or archimedean.

If $\nu$ is finite, we see that the integrand has period $2\pi\over\log(q)$ except for $\widehat{g}$. The Poisson summation formula implies
$$\eqalign{
\sum_{j\in\ZZ}\widehat{g}(j\,{2\pi i\over\log(q)} + {1\over2} + i\tau)
&= \log(q)\,\sum_{k\in\ZZ}g(q^k)\,q^{k({1\over2} + i\tau)} \cr
&= \log(q)\, \int_{K_\nu^\times}g(|t|) |t|^{{1\over2} + i\tau}\,d^*t \cr
&= \log(q)\ f(\chi\omega_\tau) }$$
so that we obtain in the end
$$F(t) = -{\log(q) \over {2\pi}}\int_{-{\pi\over\log(q)}}^{+{\pi\over\log(q)}} f(\chi\omega_\tau)\Lambda(\chi\omega_\tau, {1\over2})\,\overline{\chi\omega_\tau(t)}\, d\tau$$
Starting from
$$g(|t|) = {1 \over {2\pi i}}\int_{Re(s)=1/2} \widehat{g}(s)|t|^{-s}\, ds$$
leads in the same manner to
$$f(t) = {\log(q) \over {2\pi}}\int_{-{\pi\over\log(q)}}^{+{\pi\over\log(q)}} f(\chi\omega_\tau)\,\overline{\chi\omega_\tau(t)}\, d\tau$$
This identifies in a concrete way that component of $X_\nu$ containing $\chi(1/t)$ and proves the Theorem in the finite case.

In the archimedean case we can just identify immediately $\widehat{g}({1\over2} + i\tau)$ with $f(\chi\omega_\tau)$ as the multiplicative Haar measure $d^*t$ has been chosen to be pushed forward to $du/u$ under the norm map $t\mapsto|t|$. So the formulae in that case are
$$F(t) = -{1 \over {2\pi}}\int_{-\infty}^{+\infty} f(\chi\omega_\tau)\Lambda(\chi\omega_\tau, {1\over2})\,\overline{\chi\omega_\tau(t)}\, d\tau$$
$$f(t) = {1\over {2\pi}}\int_{-\infty}^{+\infty} f(\chi\omega_\tau)\,\overline{\chi\omega_\tau(t)}\, d\tau$$
This completes the proof of the Theorem G1.

{\bf H. More on homogeneous distributions and Gamma functions}

In the previous proof we made a detour through the complex plane. It is possible to proceed in a more economical manner, as we will show in the next chapter ({\bf I}), where we will prove the following identities of distributions on the punctured $\nu$-adics (for $0<Re(s)<1$):

$$H({\chi(x)}^{-1}|x|^{-s}) = \left({\partial\over\partial s}\log(\Gamma(\chi,s))\right) {\chi(x)}^{-1}|x|^{-s}\leqno ({\cal S}_H)$$
$$K({\chi(x)}^{-1}|x|^{-s}) = - i\,\left({\partial^2\over\partial s^2}\log(\Gamma(\chi,s))\right){\chi(x)}^{-1}|x|^{-s}\leqno ({\cal S}_K)$$

Let us examine here the consequences of these equations. The functions $H(\chi, s) = {\partial\over\partial s}\log(\Gamma(\chi,s))$ and $K(\chi, s) = - i\,{\partial^2\over\partial s^2}\log(\Gamma(\chi,s))$ have the following properties
$$H(\chi,s) = H(\chi^{-1}, 1-s)$$
$$\overline{H(\chi,s)} = H(\overline{\chi},\overline{s})$$
$$Re(s) = {1\over2}\Rightarrow H(\chi,s) = H(\overline{\chi},\overline{s})\hbox{ and is real}$$
$$K(\chi,s) = - K(\chi^{-1}, 1-s)$$
$$\overline{K(\chi,s)} = - K(\overline{\chi},\overline{s})$$
$$Re(s) = {1\over2}\Rightarrow K(\chi,s) = -K(\overline{\chi},\overline{s})\hbox{ and is real}$$
As the Inversion sends a unitary character to its complex conjugate, this implies:

{\bf Theorem H1:} The conductor (resp{.} commutator) operator commutes (resp{.} anticommutes) with the Inversion.

At a finite place the logarithmic derivative of the Gamma function of a ramified character is a constant so that the next derivative vanishes. In this manner we recover the result that $\log(|x|)$ and $\log(|y|)$ commute on the ramified part of the Hilbert space. There only remains the invariant part, which is dual to the circle of non-ramified unitary characters. The spectral functions of $H$ and $K$ are continuous hence bounded. So the commutator operator at a finite place is a bounded hermitian operator, and it remains to check that this works for the real and complex places too.

For the real (resp{.} complex) place there is a canonical way to choose a base point $\chi_\alpha$ in each component $\alpha$ of $S_r = Hom(\RR^\times, U(1))$ (resp{.} $S_c = Hom(\CC^\times, U(1))$). One requires $\chi_\alpha$ to be invariant under positive dilations, so that for the real place $r$ we have $\alpha \in\{+,-\}$ with $\chi_+(u) = 1$ and $\chi_-(u) = sgn(u)$ whereas for the complex place $c$ we have $\alpha = N, N \in \ZZ$ with $\chi_N(z) = z^N(z\overline{z})^{-N/2}$.

The associated Gamma functions are
$$\Gamma(\chi_+,s) = \gamma_+(s) = \pi^{{1\over2}-s}\,{\Gamma({s\over2})\over\Gamma({1-s\over2})}$$
$$\Gamma(\chi_-,s) = \gamma_-(s) = i\,\pi^{{1\over2}-s}\,{\Gamma({s+1\over2})\over\Gamma({2-s\over2})}$$
$$\Gamma(\chi_N,s) = \gamma_N(s) = i^{|N|}\,(2\pi)^{1-2s}\,{\Gamma({|N|\over2} + s)\over\Gamma({|N|\over2} + 1-s)}$$

In these equations $\Gamma(s)$ is Euler's Gamma function. Its logarithmic derivative $\lambda(s)$ satisfies
$$\lambda(s) = {\Gamma^\prime(s)\over\Gamma(s)} = -\gamma - {1\over s} - \sum_{j\geq1}({1\over{j+s}}-{1\over j})$$
$$\lambda(s+1) = {1\over s} + \lambda(s)$$
The spectral functions on the critical line for $H$ are
$$\eqalign{
h_+(\tau) &= H(\chi_+,{1\over2}+i\tau) = -\log(\pi) -\gamma - {1\over{1\over4}+\tau^2} + \sum_{j\geq1}\left({1\over j}-{4j+1\over(2j+{1\over2})^2 +\tau^2}\right)\cr
h_-(\tau) &= H(\chi_-,{1\over2}+i\tau) = -\log(\pi) -\gamma - {3\over{9\over4}+\tau^2} + \sum_{j\geq1}\left({1\over j}-{4j+3\over(2j+{3\over2})^2 +\tau^2}\right)\cr
h_N(\tau) &= H(\chi_N,{1\over2}+i\tau) = -2\log(2\pi) -2\gamma - {|N|+1\over({|N|+1\over2})^2+\tau^2} + 2\sum_{j\geq1}\left({1\over j}-{j+{|N|+1\over2}\over(j+{|N|+1\over2})^2 + \tau^2}\right)\cr
}$$
From this we see that they take their minimal values at $\tau=0$, are even and increase steadily to $\infty$ when $|\tau|\rightarrow\infty$ (at a logarithmic growth that can also be deduced from Stirling's Formula). These minimal values are 
$$\mu_+ = h_+(0) = -\log(\pi) + \lambda({1\over4}) = -(\log(8\pi) + \gamma) - {\pi\over2}$$
$$\mu_- = h_-(0) = -\log(\pi) + \lambda({3\over4}) = -(\log(8\pi) + \gamma) + {\pi\over2}$$
$$\mu_N = h_N(0) = -2\log(2\pi) + 2\lambda({|N|+1\over2})$$
$$\mu_{2N} = -2(\log(8\pi) + \gamma) + 2\sum_{j=1}^N{1\over j-{1\over2}}\quad (N\geq0)$$
$$\mu_{2N+1} = -2(\log(8\pi) + \gamma) + 4\log(2) + 2\sum_{j=1}^N{1\over j}\quad (N\geq0)$$

{\bf Theorem H2:} The conductor operator is bounded below.

The spectral functions of the commutator operator on the critical line are 
$$\eqalign{
k_+(\tau) &= -\sum_{j\geq0}{(4j+1)\,2\tau\over((2j+{1\over2})^2 +\tau^2)^2}\cr
k_-(\tau) &= -\sum_{j\geq0}{(4j+3)\,2\tau\over((2j+{3\over2})^2 +\tau^2)^2}\cr
k_N(\tau) &= -\sum_{j\geq0}{(2j+|N|+1)\,2\tau\over((j+{|N|+1\over2})^2 + \tau^2)^2}\cr
}$$
From this or from the Stirling Formula one finds that these functions are $O({1\over|\tau|})$ when $\tau\rightarrow\infty$, hence bounded. As $|k_{N+2}(\tau)| \leq |k_N(\tau)|$ for $N\geq0$ one gets

{\bf Theorem H3:} The commutator operator is bounded.

{\bf I. Taking derivatives and getting conductor and commutators}

The basis of the method remains the fundamental Fourier Transform identity
$${\cal F}(\chi(y)|y|^{s-1}) = \Gamma(\chi, s)\chi^{-1}(x)|x|^{-s}$$
We will always take $\chi$ to be a unitary character, then for $s$ in the critical strip this has the direct meaning that
$$\int \widetilde{\alpha}(y)\,\chi(y)|y|^{s-1} dy = \Gamma(\chi, s) \int \alpha(x)\,\chi^{-1}(x)|x|^{-s} dx\leqno{\hbox{(IF)}}$$
for any Schwartz-Bruhat function $\alpha(x)$ on $K_\nu$ and $\widetilde{\alpha}(y) = \int_{K_\nu}\varphi(x)\lambda(-xy)\,dx$ its Fourier Transform.

{\bf Lemma I1:} For $\varphi(x)$ a Schwartz-Bruhat function on $K_\nu$ its convolution $B(\varphi)$ with the Fourier Transform of $\log(|y|)$ is a continuous function which is $O({1\over|x|})$ when $|x| \rightarrow \infty$.

{\bf Proof:} As $B(\varphi)$ is the Fourier transform of an $L^1$ function it is continuous. The formula for the Fourier transform of $\log(|y|)$ was given above. For a finite place it implies that $B(\varphi)(x)$ for $|x|$ large enough is identically ${C\over|x|}$ for some constant $C(\varphi)$. At the real place we write it as
$$-B(\varphi)(x) = \int_{|t| \leq 1}(\varphi(x-t) - \varphi(x))\, {dt \over 2|t|} + \int_{|t| > 1}\varphi(x-t)\, {dt \over 2|t|}+\ (\log(2\pi) + \gamma)\cdot\varphi(x)$$
The first integral is bounded by $\sup_{|t|\geq|x|-1}(|\varphi^\prime(t)|)$. For the second integral one first looks at the contribution of the region of integration $|t| > |x|/2$ which clearly gives $O({1\over|x|})$, and then the remaining part is bounded by $|x|\sup_{|t|\geq|x|/2}(|\varphi(t)|)$ hence with rapid decrease as $|x|\rightarrow\infty$. The complex case is similar and left to the reader (recall though that $|x|$ then means $z\overline{z}$ for $x=z\in\CC$).

{\bf Lemma I2:} The identity for $s$ in the critical strip
$$\int \psi(y)\,\chi(y)|y|^{s-1} dy = \Gamma(\chi, s) \int \varphi(x)\,\chi^{-1}(x)|x|^{-s} dx$$
where both $\varphi(x)$ and $\psi(y)$ are supposed to be measurable locally integrable functions, tempered as distributions, and $\psi(y)$ the Fourier transform of $\varphi(x)$ as a distribution, is valid as soon as both integrals make sense as Lebesgue Integrals (that is are absolutely convergent).

{\bf Proof:} One just needs to check that the double integral and change of variables trick of Tate's Thesis works. Let $\alpha(x)$ be a Schwartz-Bruhat function with Fourier Transform $\widetilde{\alpha}(y)$. Starting with
$$\int \psi(y)\,\chi(y)|y|^{s-1} dy \cdot \int \alpha(x)\,\chi^{-1}(x)|x|^{-s} dx$$
we convert it by Fubini's Theorem to a double integral and then apply the change of variables $y = uv, x =v$ to get
$$\int\int\psi(uv)\alpha(v)\chi(u)|u|^{s-1}\,dudv$$
Then one evaluates for $u\neq0$
$$\int\psi(uv)\alpha(v)dv = {1\over|u|}\int\varphi({w\over u})\widetilde{\alpha}(w)dw = \int\varphi(v)\widetilde{\alpha}(uv)\,dv$$
so that the double integral becomes
$$\int\int\varphi(v)\widetilde{\alpha}(uv)\chi(u)|u|^{s-1}\,dudv$$
which is similarly seen (working in reverse) to be
$$\int \widetilde{\alpha}(y)\,\chi(y)|y|^{s-1} dy \cdot \int \varphi(x)\,\chi^{-1}(x)|x|^{-s} dx$$
We have only used the definition of the Fourier Transform of a distribution and the fact that all integrals considered were absolutely convergent, so that the manipulations are allowed by Fubini's Theorem.

Let's start from the fundamental identity (IF). We are allowed by Lemmas I1 and I2 to do two things: to take its derivative, and on the other hand to apply it directly to the Fourier pair $\beta(x) = B(\alpha)(x), \widetilde{\beta}(y) = \log(|y|)\widetilde{\alpha}(y)$. The first operation (with $H(\chi,s)$ defined to be the logarithmic derivative of $\Gamma(\chi,s)$) gives
$$H(\chi,s)\int_{K_\nu}\alpha(x)\chi^{-1}(x)|x|^{-s}\,dx = \int_{K_\nu}\alpha(x)\log(|x|)\chi^{-1}(x)|x|^{-s}\,dx + \int_{K_\nu}\widetilde{\alpha}(y)\log(|y|){\chi(y)|y|^{s-1}\over\Gamma(\chi,s)}\,dy$$
and the second gives:
$$\int_{K_\nu}\widetilde{\alpha}(y)\log(|y|){\chi(y)|y|^{s-1}\over\Gamma(\chi,s)}\,dy
= \int_{K_\nu}B(\alpha)(x)\chi^{-1}(x)|x|^{-s}\,dx$$
so that we end up with
$$H(\chi,s)\int_{K_\nu}\alpha(x)\chi^{-1}(x)|x|^{-s}\,dx = \int_{K_\nu}H(\alpha)(x)\chi^{-1}(x)|x|^{-s}\,dx\leqno{\hbox{(IH)}}$$
Furthermore by Lemma I1, the following integral
$$\int_{K_\nu}H(\alpha_1)(x)\alpha_2(x)\,dx$$
for two test-functions $\alpha_1$ and $\alpha_2$ is absolutely convergent and easily seen to be
$$\int_{K_\nu}\alpha_1(x)H(\alpha_2)(x)\,dx$$
so that the Identity (IF) can be interpreted as an identity of distributions
$$H({\chi(x)}^{-1}|x|^{-s}) = H(\chi, s)\ {\chi(x)}^{-1}|x|^{-s}
 = \left({\partial\over\partial s}\log(\Gamma(\chi,s))\right)\ {\chi(x)}^{-1}|x|^{-s}\leqno ({\cal S}_H)$$

To proceed further we now restrict $\alpha(x)$ to vanish identically in a neighborhood of the origin so that $A(\alpha)(x) = \log(|x|)\alpha(x)$ is again a Schwartz function. Again on one hand we compute the derivative of (IH)
$$H^\prime(\chi,s)\int_{K_\nu}\alpha(x)\chi^{-1}(x)|x|^{-s}\,dx$$
$$= H(\chi,s)\int_{K_\nu}\alpha(x)\log(|x|)\chi^{-1}(x)|x|^{-s}\,dx -
\int_{K_\nu}H(\alpha)(x)\log(|x|)\chi^{-1}(x)|x|^{-s}\,dx$$
and on the other hand we apply it directly to $A(\alpha)(x)$:
$$H(\chi,s)\int_{K_\nu}\alpha(x)\log(|x|)\chi^{-1}(x)|x|^{-s}\,dx
= \int_{K_\nu}H(A(\alpha))(x)\chi^{-1}(x)|x|^{-s}\,dx$$
so that we end up with
$$H^\prime(\chi,s)\int_{K_\nu}\alpha(x)\chi^{-1}(x)|x|^{-s}\,dx = \int_{K_\nu}((HA-AH)(\alpha))(x)\chi^{-1}(x)|x|^{-s}\,dx$$
hence with
$$iH^\prime(\chi,s)\int_{K_\nu}\alpha(x)\chi^{-1}(x)|x|^{-s}\,dx = \int_{K_\nu}(K(\alpha))(x)|x|^{-s}\,dx\leqno{\hbox{(IK)}}$$
For $\alpha_1$ and $\alpha_2$ vanishing in a neighborhood of the origin one checks easily
$$\int_{K_\nu}K(\alpha_1)(x)\alpha_2(x)\,dx = - \int_{K_\nu}\alpha_1(x)K(\alpha_2)(x)\,dx$$
so that we can reinterpret (IK) as an identity of distributions on the {\it punctured} $\nu$-adics
$$K({\chi(x)}^{-1}|x|^{-s}) = K(\chi, s)\ {\chi(x)}^{-1}|x|^{-s}
= \left(- i\,{\partial^2\over\partial s^2}\log(\Gamma(\chi,s))\right)\ {\chi(x)}^{-1}|x|^{-s}\leqno ({\cal S}_K)$$

Let's now proceed inductively with higher derivatives and hence deduce the spectral analysis of the higher commutators $K_1 = iK = [A,H], K_{N+1} = [A, K_N]$ as an identity of distributions on the punctured $\nu$-adics
$$K_N({\chi(x)}^{-1}|x|^{-s}) = \left({\partial^{N+1}\over\partial s^{N+1}}\log(\Gamma(\chi,s))\right)\ {\chi(x)}^{-1}|x|^{-s}\leqno ({\cal S}_N)$$
from which the Hilbertian version follows by restriction to the critical line (note that these operators are indeed completely dilation invariant as any commutator built with the identity vanishes).

We apply the operators $K_N$ only to Schwartz functions vanishing identically in a neighborhood of the origin. The explicit formula
$$K_N = \sum_{N\geq j\geq0}(-1)^{N-j}{N\choose j}A^jHA^{N-j}$$
shows that the resulting function is (according to Lemma I1) continuous except perhaps at the origin where it is $O(\log(|x|)^N)$ while at infinity it is $O({\log(|x|)^N\over|x|})$ so that Lemma I2 allows to make sense of $K_N({\chi(x)}^{-1}|x|^{-s})$ as a distribution on the punctured $\nu-$adics and to generalize the method used for $H$ and $K$.

Assuming inductively the validity of:
$$\left({\partial^{N+1}\over\partial s^{N+1}}\log(\Gamma(\chi,s))\right)\int_{K_\nu}\alpha(x)\chi^{-1}(x)|x|^{-s}\,dx = (-1)^N\int_{K_\nu}(K_N(\alpha))(x)\chi^{-1}(x)|x|^{-s}\,dx\leqno{\hbox{(I(N))}}$$
we first compute its derivative
$$\left({\partial^{N+2}\over\partial s^{N+2}}\log(\Gamma(\chi,s))\right)\int_{K_\nu}\alpha(x)\chi^{-1}(x)|x|^{-s}\,dx =$$
$$\left({\partial^{N+1}\over\partial s^{N+1}}\log(\Gamma(\chi,s))\right)\int_{K_\nu}\alpha(x)\log(|x|)\chi^{-1}(x)|x|^{-s}\,dx - (-1)^N\int_{K_\nu}(K_N(\alpha))(x)\log(|x|)|x|^{-s}\,dx$$
and then apply (I(N)) to $A(\alpha)$ so that
$$\left({\partial^{N+2}\over\partial s^{N+2}}\log(\Gamma(\chi,s))\right)\int_{K_\nu}\alpha(x)\chi^{-1}(x)|x|^{-s}\,dx
= (-1)^{N+1}\int_{K_\nu}([A,K_N](\alpha))(x)|x|^{-s}\,dx$$
which is (I(N+1)).

As for $\alpha_1$ and $\alpha_2$ two test-functions with support away from the origin 
$$\int_{K_\nu}K_N(\alpha_1)(x)\alpha_2(x)\,dx = (-1)^N \int_{K_\nu}\alpha_1(x)K_N(\alpha_2)(x)\,dx$$
this gives the identity of distributions $({\cal S}_N)$ on the punctured $\nu-$adics.

As we did before for $K$ from the partial fraction expansion of Euler's Gamma function we see that from the Hilbert Space point of view all $K_N$'s are bounded.

{\bf J. More on the Hilbert Space story}

Let $\Delta$ be the domain of functions $f(t)$ such that the unitary dual $f(\chi)$ in $L^2(X_\nu, d\mu)$ is of the Schwartz class in each component (which means smooth for a finite place as the components are then circles) and non-zero in only finitely many components. In the additive picture this translates into\ : if the place is finite, with $q$ the cardinality of the residue field, $\Delta$ consists of functions $\varphi(x) = |x|^{-1/2}F({\log(|x|)\over\log(q)})$ where the $F(n), n\in\ZZ$ are the Fourier coefficients of a smooth function on the circle (that is decrease  faster than any polynomial in $|n|^{-1}$). We also allow twists with multiplicative unitary characters, and then any finite linear combination of such building blocks. At the real or complex place this means $\varphi(x) = |x|^{-1/2}F(\log(|x|))$ or a twist with a unitary character, and finite linear combinations, where $F$ is a Schwartz function on $\RR$. 

{\bf Theorem J1:} The dense subdomain $\Delta$ of $L^2(K_\nu, dx)$ is stable under the Fourier Transform, $A$, $B$, $H$, $K$, and all higher $K_N$'s. All commutators on this domain of $A$ and $B$ boil down to the already constructed $K_N$'s.

{\bf Proof:} Clearly this $\Delta$ is stable under multiplication with $\log(|x|)$. And it is also stable under additive convolution with the Fourier Transform of $\log(|y|)$ ! Indeed it is enough to show that it is stable under the (additive) Fourier Transform, but the Fourier Transform was seen to act in $L^2(X_\nu, d\mu)$ through a composition of Inversion ($\chi \mapsto \overline{\chi}$) and multiplication with the Gamma function (seen on the critical line). For a finite place the Gamma function is smooth on each component so we are done, for an archimedean place all derivatives of $\Gamma(\chi,{1\over2}+i\tau)$ have a growth bounded by a polynomial in $\log(|\tau|)$. This is seen inductively from $\Gamma^\prime(\chi,s) = H(\chi,s)\Gamma(\chi,s)$ as the $H$ term has logarithmic growth and all its derivatives are bounded (as seen before, this is easily deduced from the partial fraction expansion of Euler's Gamma function), whereas $\Gamma(\chi,s)$ itself is (of course) of absolute value $1$. So multiplying with this leaves the Schwartz space invariant. It only remains to show that no new commutator can be built. But obviously $[H,K_N] = 0$ and $[K_N, K_M] = 0$ from the spectral picture. So any commutator built with $A$'s and $B$'s can be expressed in terms of $A$'s and only one $B$, that is it is one of the $K_N$'s.

{\bf Theorem J2:} The domain ${\cal S}_0$ is a domain of self-adjointness for the symmetric operator $H$.

{\bf Proof:} The dense domain $\Delta$ is stable under $H$ and in fact $H \pm i$ is boundedly invertible on it, so it is a domain of self-adjointness for $H$. To conclude we need to prove that for any $f$ in $\Delta$ we can find $f_0$ in ${\cal S}_0$ making $||f - (H + i)(f_0)||$ arbitrarily small (and the same with $H - i$). This will follow if for any $f$ in $\Delta$ we can find $f_0$ in ${\cal S}_0$ making $||(H + i)(f - f_0)||$ arbitrarily small. Working in the spectral picture, we see that this is clearly possible if the place is non-archimedean as the components of $X_\nu$ are then compact. The only slight difficulty is at an archimedean place, for concreteness let us consider the case of the real place and of an even $f$, the other cases being disposed of similarly. We first bound $|h_+(\tau) + i|^2$ with $O(1 + \tau^2)$. Taking Fourier Transforms then reduces the problem to the question of approximating in the sense of $\int_\RR \left(|\psi(a)|^2 + |\psi'(a)|^2\right)\, da$ an arbitrary Schwartz function $\psi(a)$ with a compactly supported one. One just takes $\psi_0(a) = \rho({a\over\Lambda})\psi(a)$ with $\rho$ smooth, $\rho(a) = 1$ for $|a|\leq 1$, $\rho(a) = 0$ for $|a| > 2$, and $\Lambda \rightarrow \infty$.

If the goal is just to understand $H$ and $K$ in the Hilbert space, the use of the domain $\Delta$ allows for a rather expeditious method. For $\chi$ a fixed unitary character we parametrize the corresponding component of $X_\nu$ through the association of the character $\chi(u)|u|^{i\tau}$ to the real number $\tau$, taken modulo ${2\pi i \over\log(q)}$ if the place is finite. An element $f(t)$ of $L^2(K_\nu^\times, d^*t)$ has a multiplicative Fourier Transform $f(\chi,\tau) = \int_{K_\nu^\times}f(t)\chi(t)|t|^{i\tau}\,d^*t$. We can then take the derivative with respect to $\tau$ and end up with:

{\bf Theorem J3:} In this picture $A$ is ${1\over i}{\partial\over\partial\tau}$

We know how the additive Fourier Transform acts:
$$f(\chi,\tau)\mapsto \Gamma(\chi,{1\over2}+i\tau)f(\overline{\chi},-\tau)$$
The inverse Fourier transform acts as
$$f(\chi,\tau)\mapsto \Gamma(\chi,{1\over2}+i\tau)\chi(-1)f(\overline{\chi},-\tau)$$
so $B = {\cal F}A{\cal F}^{-1}$ maps $f(\chi,\tau)$ to
$$\Gamma(\chi,{1\over2}+i\tau)\Gamma^\prime(\overline{\chi},{1\over2}-i\tau)\chi(-1)f(\chi,\tau)
+\Gamma(\chi,{1\over2}+i\tau)\Gamma(\overline{\chi},{1\over2}-i\tau)\chi(-1){1\over i}(-1){\partial\over\partial\tau}f(\chi,\tau)$$
This can be simplified using $\Gamma(\chi,s)\Gamma(\chi^{-1},1-s)=\chi(-1)$ to 
$$H(\chi,{1\over2}+i\tau)f(\chi,\tau) - {1\over i}{\partial\over\partial\tau}f(\chi,\tau)$$
So the spectral function of $H = A+B$ is as was previously found (as $H(\chi,{1\over2}+i\tau) = H(\overline{\chi},{1\over2}-i\tau)$), and all commutators follow from Theorem J3 as well.

{\bf K. Notes}

{\bf A.} The material in this paper is covered in my electronic manuscripts [Bu98a], [Bu98b], [Bu98c], where some additional comments will be found. A related manuscript is [Bu98d].

{\bf B.} The original impetus was the similarity between the local terms of the Explicit Formula and the $n$-point scattering amplitudes of the open string. On the subject of the open string over the $p$-adic field, I refer the reader to the beautiful paper of Zabrodin [Za].

{\bf C.} It is interesting to note that Weil mentions the Fourier Transform of $\log(|t|)$ in [We52], in the context of applying Parseval's formula to the integral on the critical line involving the logarithmic derivative of Euler's Gamma function. He says about it that it is a distribution that has been computed by Laurent Schwartz, but as far as I know, does not mention the even closer connection with his final result.

{\bf D.} The problem of treating all local terms in the most uniform manner possible has also been addressed by Connes [Co98] (among many other things). His method is very geometric and leads to Weil's terms in their original form. One notices a confluence of ideas on the importance of looking at the additive Fourier Transform in multiplicative terms. In [Bu99] I apply the Transformation Theory to the calculation of a Hilbert trace first evaluated by Connes [Co98].

{\bf E.} In [BeKe98a, BeKe98b], Berry and Keating have applied techniques from mathematical physics to the operator $H = xp$ and discovered a connection with the smooth part of the density of zeros of the Riemann Zeta function. They emphasize the necessity for this to identify $x$ with $-x$, $p$ with $-p$ and the role played by the Inversion (called by them ``quantum exchange''). It is quite satisfying to see that the conductor operator $\log(|x|) + \log(|p|)$ has these symmetries. But the way it exactly relates to the work of Berry and Keating and their further stimulating suggestions remains to be clarified. In [Bu99] I show that at a non-archimedean place it has an interpretation in terms of scattering.

{\bf F.} The boundedness of the commutator operator should be related to the theory of band-and-time limited functions [SP61], mentioned in Connes's paper. It would also be interesting to know the exact operator norm of the commutator operator at the archimedean places.

{\bf G.} The set of global additive characters on the adeles such that the diagonally embedded number field $K$ is its own annihilator is a torsor under the multiplicative group $K^\times$. The local terms of the Explicit Formula are not invariant under such a change (although their sum of course is). So $K^\times$ is a non-trivial symmetry of the Explicit Formula. Another symmetry is the Inversion. It is to be noted that Goldfeld has proposed an approach to the Explicit Formula centered around the group generated by $K^\times$ and the Inversion [Go94].

\par\vfill\eject

{{\bf REFERENCES}\par
\baselineskip = 14 pt
\parskip = 6 pt
\font\smallRoman = cmr10
\smallRoman
\font\smallBold = cmbx10
\font\smallSlanted = cmsl10

{\smallBold [Ba81] K. Barner}, {\smallSlanted ``On A.Weil's explicit formula''}, Journal f\"ur Mathematik Band 323, 139-152, (1981).\par

{\smallBold [BeKe98a] M.V. Berry, J.P. Keating}, {\smallSlanted ``H = xp and the Riemann zeros''}, in {\smallSlanted ``Supersymmetry and trace formulae: chaos and disorder''}, eds. J.P. Keating, D.E. Khmelnitskii, and I.V. Lerner, Plenum, New York (1998)

{\smallBold [BeKe98b] M.V. Berry, J.P. Keating}, {\smallSlanted ``The Riemann zeros and eigenvalue asymptotics''}, HP Laboratories Technical Report HPL-BRIMS-98-26 (December 1998)

{\smallBold [Bu98a] J.F. Burnol}, {\smallSlanted ``The Explicit Formula and a Propagator''}, math/9809119 (v2 November 1998)\par

{\smallBold [Bu98b] J.F. Burnol}, {\smallSlanted ``Spectral analysis of the local conductor operator''}, math/9811040 (November 1998)\par

{\smallBold [Bu98c] J.F. Burnol}, {\smallSlanted ``Spectral analysis of the local commutator operators''}, math/9812012 (December 1998)\par

{\smallBold [Bu98d] J.F. Burnol}, {\smallSlanted ``The Explicit Formula in simple terms''}, math/9810169 (v2 November 1998).\par

{\smallBold [Bu99] J.F. Burnol}, {\smallSlanted ``Scattering on the $p-$adic field and a trace formula''}, math/9901051 (January 1999).\par

{\smallBold [Co98] A. Connes}, {\smallSlanted ``Trace formula in non-commutative Geometry and the zeros of the Riemann zeta function''}, math/9811068 (October 1998).\par

{\smallBold [GeGrPi69] I. M. Gel'fand, M. I. Graev, I. I. Piateskii-Shapiro},{\smallSlanted ``Representation Theory and automorphic functions''}, Philadelphia, Saunders (1969).\par

{\smallBold [Go94] D. Goldfeld}, {\smallSlanted ``A spectral interpretation of Weil's explicit formula''}, Lect. Notes Math 1593, p 135-152, Springer Verlag (1994).\par

{\smallBold [Ha90] S. Haran},{\smallSlanted ``Riesz potentials and explicit sums in arithmetic''}, Invent. Math.  101, 697-703 (1990).\par

{\smallBold [Iw52] K. Iwasawa}, {\smallSlanted ``Letter to J. Dieudonn\'e'' (1952)}, in {\smallSlanted ``Zeta Functions in Geometry''}, Adv.Stud.Pure Math 21 ed. Kurokawa, Sunada, pub. Kinokuniya, (1992).\par

{\smallBold [SP61] D. Slepian, H. Pollak},{\smallSlanted ``Prolate spheroidal wave functions, Fourier analysis and uncertainty I''}, Bell Syst. Tech. J. 40, (1961).\par

{\smallBold [Ta50] J. Tate}, {\smallSlanted ``Fourier Analysis in Number Fields and Hecke's Zeta Function''}, Princeton 1950, reprinted in {\smallSlanted ``Algebraic Number Theory''}, ed. J.W.S. Cassels and A. Fr\"ohlich, Academic Press, (1967).\par

{\smallBold [Vl88] V. S. Vladimirov},{\smallSlanted ``Generalized functions over the fields of p-adic numbers''}, Russian Math. Surveys 43:5, 19-64 (1988).\par

{\smallBold [We52] A. Weil},{\smallSlanted ``Sur les ``formules explicites'' de la th\'eorie des nombres premiers''}, Comm. Lund (vol d\'edi\'e \`a Marcel Riesz), (1952).\par

{\smallBold [We66] A. Weil},{\smallSlanted ``Fonctions z\^etas et distributions''}, S\'eminaire Bourbaki n${}^{\smallRoman o}$ 312, (1966).\par

{\smallBold [We74] A. Weil},{\smallSlanted ``Basic Number Theory''}, Springer Verlag, $2^{\smallRoman nd}$ edition, New York(1974).\par

{\smallBold [Za89] A. V. Zabrodin},{\smallSlanted ``Non-Archimedean Strings and Bruhat-Tits Trees''}, Commun. Math. Phys. 123, 463-483 (1989).

\vfill
\centerline{Jean-Fran\c{c}ois Burnol, 62 rue Albert Joly, F-78000 Versailles, France}
\centerline{February 1999}
}
\eject
\bye